\documentclass[12pt]{article}
\usepackage[reqno]{amsmath}
\usepackage{amssymb, theorem, enumerate, bm}
\usepackage{url}
\usepackage{hyperref}
\hypersetup{colorlinks=false}

\newcommand\cF{{\mathcal F}}
\newcommand\cG{{\mathcal G}}
\newcommand\cH{{\mathcal H}}

\newcommand\cS{{\mathcal S}}

\makeatletter
\g@addto@macro\bfseries{\boldmath}
\makeatother

\expandafter\ifx\csname pdfoptionalwaysusepdfpagebox\endcsname\relax\else
\fi

\theoremstyle{plain}
{\theorembodyfont{\slshape}
\newtheorem{theorem}{Theorem}[section]
\newtheorem{lemma}[theorem]{Lemma}
}
\theorembodyfont{\rmfamily}

\newtheorem{conjecture}[theorem]{Conjecture}

\newcommand\lref[1]{Lemma~\ref{lem:#1}}
\newcommand\tref[1]{Theorem~\ref{thm:#1}}
\newcommand\cref[1]{Corollary~\ref{cor:#1}}

\newcommand\conjref[1]{Conjecture~\ref{conj:#1}}

\def\sqr#1#2{{\vbox{\hrule height.#2pt
    \hbox{\vrule width.#2pt height#1pt \kern#1pt
        \vrule width.#2pt}\hrule height.#2pt}}}
\def\eqed{\sqr53}
\def\qed{%
    \ifmmode\eqno\eqed
    \else\nobreak\ \hfill\eqed\medbreak\fi}

\newcommand{\upperRomannumeral}[1]{\uppercase\expandafter{\romannumeral#1}}

\title{A New Quadratic Bound for the Manickam--Mikl\'{o}s--Singhi Conjecture}
\author{Ameera Chowdhury \thanks{Department of Mathematical Sciences, Carnegie Mellon University, Pittsburgh, PA, 15213, USA. E-mail: {\tt ameera@math.cmu.edu.} Research supported by NSF grant DMS-1203982.} \and Ghassan Sarkis \thanks{Department of Mathematics, Pomona College, Claremont, CA, 91711, USA. E-mail: {\tt Ghassan.Sarkis@pomona.edu, SShahriari@pomona.edu}} \and Shahriar Shahriari \footnotemark[2]}

\begin{document}
\maketitle

\begin{abstract}
More than twenty-five years ago, Manickam, Mikl\'{o}s, and Singhi conjectured that for positive integers $n,k$ with $n \geq 4k$, every set of $n$ real numbers with nonnegative sum has at least $\binom{n-1}{k-1}$ $k$-element subsets whose sum is also nonnegative. We verify this conjecture when $n \geq 8k^{2}$, which simultaneously improves and simplifies a bound of Alon, Huang, and Sudakov and also a bound of Pokrovskiy when $k < 10^{45}$.
\end{abstract}

\section{Introduction}
\label{sec:intro}

Manickam, Mikl\'{o}s, and Singhi \cite{mm, ms} conjectured in 1988 that
\begin{conjecture}
\label{conj:MMS}
For positive integers $n,k \in \mathbb{Z}^{+}$ with $n \geq 4k$, every set of $n$ real numbers with nonnegative sum has at least $\binom{n-1}{k-1}$ $k$-element subsets whose sum is also nonnegative.
\end{conjecture}
\conjref{MMS} was motivated by studies of the first distribution invariant in certain association schemes, and may also be considered an analogue of the Erd\H{o}s--Ko--Rado theorem \cite{ekr}. The Erd\H{o}s--Ko--Rado theorem states that if $n > 2k$, then any family of $k$-element subsets of an $n$-element set with the property that any two subsets have nonempty intersection has size at most $\binom{n-1}{k-1}$; moreover the unique extremal family is a star, the family of $k$-element subsets containing a fixed element.

\conjref{MMS} is similar to the Erd\H{o}s--Ko--Rado theorem, not only in the appearance of the binomial coefficient $\binom{n-1}{k-1}$, but also because the family of $k$-element subsets with nonnegative sum attains this lower bound and forms a star when one of the $n$ real numbers equals $n-1$ and the remaining $n-1$ numbers equal $-1$. As in the Erd\H{o}s--Ko--Rado theorem, $n$ must be large enough with respect to $k$, otherwise there exist $n$ real numbers with nonnegative sum and fewer than $\binom{n-1}{k-1}$ $k$-element subsets with nonnegative sum. Such examples can be easily constructed when $n=3k+r$ and $1 \leq r \leq k/7$. Although \conjref{MMS} and the Erd\H{o}s--Ko--Rado theorem share the same bound and extremal example, there is no obvious way to translate one question into the other.

\conjref{MMS} has attracted a lot of attention due to its connections with the Erd\H{o}s--Ko--Rado theorem \cite{AAH, ahs, aydblin, bhatt, bier, biermanickam, blini, chowdhury, CSS, frankl, HS, HuangSudakov, manickam, mm, ms, MC, Pokrovskiy, Tyomkyn}, but still remains open. For more than two decades, \conjref{MMS} was known to hold only when $k|n$ \cite{ms} or when $n$ is at least an exponential function of $k$ \cite{bhatt, biermanickam, mm, Tyomkyn}. In their recent breakthrough paper, Alon, Huang, and Sudakov \cite{ahs} obtained the first polynomial bound $n \geq \min\{33k^{2}, 2k^{3}\}$ on \conjref{MMS}. Later, Aydinian and Blinovsky \cite{aydblin} and Frankl \cite{frankl} gave different proofs of \conjref{MMS} for a cubic range. Recently, a linear bound $n > 10^{46}k$ has been obtained by Pokrovskiy \cite{Pokrovskiy}. Finally, there are also several works that verify \conjref{MMS} for small $k$ \cite{chowdhury, HS, manickam, MC}.

The main result of this paper verifies \conjref{MMS} when $n \geq 8k^{2}$. In particular, \tref{quadratic} simultaneously improves and simplifies the bound $n \geq \min\{33k^{2}, 2k^{3} \}$ of Alon, Huang, and Sudakov \cite{ahs} and also the bound $n \geq 10^{46}k$ of Pokrovskiy \cite{Pokrovskiy} when $k < 10^{45}$. Note that there is no loss of generality in assuming that the $n$ real numbers in \conjref{MMS} are listed in decreasing order and sum to zero.

\begin{theorem}
\label{thm:quadratic}
Let $X = \{x_{1}, \ldots, x_{n} \} \subset \mathbb{R}$ be a set of $n$ real numbers whose sum is zero, and assume $x_{i} \geq x_{j}$ if $i \leq j$. If $n \geq 8k^{2}$, then at least $\binom{n-1}{k-1}$ $k$-element subsets of $X$ have nonnegative sum. Moreover, if equality holds, the family of $k$-element subsets with nonnegative sum is a star on $x_{1}$, $\left \{S \in \binom{X}{k} : x_{1} \in S \right \}$.
\end{theorem}

The proof of \tref{quadratic} is similar to that of Theorem 1.3 in \cite{CSS}, where we tackle the Manickam--Mikl\'{o}s--Singhi conjectures for sets and vector spaces simultaneously. For the reader's convenience, we present the calculations for the case of sets in full detail in this unpublished manuscript.

\section{Bose-Mesner Matrices}
\label{sec:bosemesner}

We will need a lemma involving inclusion matrices $W_{jk}$ and Kneser matrices $\overline{W}_{jk}$. Let $W_{jk}$ (respectively $\overline{W}_{jk}$) denote the $\binom{n}{j} \times \binom{n}{k}$ matrix whose rows are indexed by the $j$-element subsets of $X$, whose columns are indexed by the $k$-element subsets of $X$, and where the entry in row $Y$ and column $S$ is $1$ if $Y \subset S$ (respectively if $Y \cap S = \emptyset$) and is 0 otherwise.

Define $\vec{x} = (x_{1}, \ldots, x_{n}) \in \mathbb{R}^{[n]}$ to be a vector that lists the $n$ real numbers in \tref{quadratic}. Without loss of generality, we may assume that $\vec{x} \neq \vec{0}$, and observe that $\vec{x}$ is orthogonal to $\vec{1}$. We will show that $W_{1k}^{T} \vec{x}$ has at least $\binom{n-1}{k-1}$ nonnegative entries when $n \geq 8k^{2}$.

An important observation by Wilson \cite{WilsonExact} is that $W_{1k}^{T} \vec{x}$ is an eigenvector of the Bose-Mesner matrix
\begin{equation}
\label{bosemesner}
B_{j} = \overline{W}_{jk}^{T}W_{jk}
\end{equation}
for $0 \leq j \leq k$ with eigenvalue $- \binom{k-1}{j-1} \binom{n-j-1}{k-1}$. We include a proof for completeness.

\begin{lemma}[Wilson, \cite{WilsonExact}]
\label{lem:eigenvecset}
For $0 \leq j \leq k$, we have $W_{1k}^{T} \vec{x}$ is an eigenvector of the Bose-Mesner matrix $B_{j}$ with eigenvalue
\begin{equation}
\label{seteigenval}
-\binom{k-1}{j-1} \binom{n-j-1}{k-1}.
\end{equation}
\end{lemma}

\noindent \textbf{Proof.} Since the columns of $W_{1k}^{T}$ are linearly independent \cite{GJ, Kantor, WilsonNecessary} and $\vec{x} \neq \vec{0}$, we have that $W_{1k}^{T} \vec{x} \neq \vec{0}$. For $S \in \binom{X}{j}$ and $T \in \binom{X}{1}$, observe that $W_{jk}W_{1k}^{T}(S,T)$ counts the number of $k$-element subsets of $X$ that contain $S \cup T$. Hence,
\begin{equation}
\label{SunionT}
W_{jk}W_{1k}^{T}(S,T) = \begin{cases}
                        \binom{n-j}{k-j} & \text{if $T \subset S$} \\
                        \binom{n-j-1}{k-j-1} & \text{if $T \not \subset S$.} \\
                        \end{cases}
\end{equation}
For the remainder of this proof, $J$ is a matrix all of whose $n$ columns are $\vec{1}$. Since $J \vec{x} = \vec{0}$, we have
\begin{equation}
\label{j1}
W_{jk}W_{1k}^{T}\vec{x} = \left( \binom{n-j-1}{k-j} W_{1j}^{T} + \binom{n-j-1}{k-j-1}J \right)\vec{x} = \binom{n-j-1}{k-j}W_{1j}^{T} \vec{x}.
\end{equation}
For $A \in \binom{X}{k}$ and $B \in \binom{X}{1}$, observe that $\overline{W}_{jk}^{T} W_{1j}^{T}(A,B)$ counts the number of $j$-element subsets of $X$ that are disjoint from $A$ and that contain $B$. Hence,
\begin{equation}
\label{kk}
\overline{W}_{jk}^{T}W_{1j}^{T} \vec{x} = \binom{n-k-1}{j-1} \overline{W}_{1k}^{T} \vec{x} = \binom{n-k-1}{j-1}(J - W_{1k}^{T}) \vec{x} = -\binom{n-k-1}{j-1}W_{1k}^{T} \vec{x},
\end{equation}
since $J \vec{x} = \vec{0}$.
Multiplying \eqref{j1} on the left by $\overline{W}_{jk}^{T}$ and applying \eqref{kk} yields
\begin{equation}
\label{eigenval}
B_{j}W_{1k}^{T} \vec{x} = - \binom{n-j-1}{k-j} \binom{n-k-1}{j-1}W_{1k}^{T} \vec{x} = - \binom{k-1}{j-1} \binom{n-j-1}{k-1}W_{1k}^{T} \vec{x},
\end{equation}
which proves that $W_{1k}^{T} \vec{x}$ is an eigenvector of the Bose-Mesner matrix $B_{j}$ with eigenvalue $-\binom{k-1}{j-1} \binom{n-j-1}{k-1}$. \qed

We will use \lref{eigenvecset} to obtain lower bounds on the number of nonnegative $k$-element subsets that intersect $\{x_{1}, \ldots, x_{k} \}$ and that contain $x_{1}$ respectively.

\section{Bounds from Eigenvalues}
\label{sec:linear}

The main result of this section is \lref{lotson1sets}, which shows that if $n \geq Ck^{2}$ and there are at most $\binom{n-1}{k-1}$ nonnegative $k$-element subsets of $X$, then at least $(1 - \frac{6}{C}) \binom{n-1}{k-1}$ $k$-element subsets on $x_{1}$ have nonnegative sum.

Henceforth, we write $W_{1k}^{T} \vec{x} = \vec{b}$, and index the entries of $\vec{b}$ with subsets $S \in \binom{X}{k}$.
Let $A = \{x_{1}, \ldots, x_{k}\}$ and note that $A$ is the $k$-element subset of $X$ with largest sum. Using \lref{eigenvecset}, we give a lower bound on the number of nonnegative $k$-element subsets of $X$ that intersect $A$.

\begin{lemma}
\label{lem:highestweightintersectionset}
There are greater than $\binom{n-k-1}{k-1}$ nonnegative $k$-element subsets of $X$ that intersect $A = \{x_{1}, \ldots, x_{k} \}$.
\end{lemma}

\noindent \textbf{Proof.}  Observe that $b_{A}$ is a largest entry of $\vec{b}$ and that $b_{A} > 0$ since $\vec{b} \neq \vec{0}$ and $\vec{b}$ is orthogonal to $\vec{1}$.

For $S, T \in \binom{X}{k}$, observe that $B_{j}(S,T)$ counts the number of $j$-element subsets of $X$ that lie in $T$ and are disjoint from $S$. Hence,
\begin{equation}
\label{bjentries}
B_{j}(S,T) = \binom{k-|S \cap T|}{j}.
\end{equation}

By \lref{eigenvecset}, the dot product of the row of $B_{k}$ corresponding to $A$ and $\vec{b}$ equals $-\binom{n-k-1}{k-1} b_{A}$. Hence, \eqref{bjentries} with $j=k$ yields
\begin{equation}
\label{sumonempty}
\sum_{S \cap A = \emptyset} b_{S} = -\binom{n-k-1}{k-1} b_{A}.
\end{equation}
Since $\vec{b}$ is orthogonal to $\vec{1}$, we see that
\begin{equation}
\label{orthogonal}
\sum_{S \cap A \neq \emptyset} b_{S} = \binom{n-k-1}{k-1} b_{A}.
\end{equation}
Since $b_{A}$ is a largest entry of $\vec{b}$, there are greater than $\binom{n-k-1}{k-1}$ nonnegative $k$-element subsets of $X$ that intersect $A$. \qed

Recall that $A = \{x_{1}, \ldots, x_{k} \}$ is the $k$-element subset of $X$ with largest sum. Let $C = \{x_{1}, x_{k+1}, \ldots, x_{2k-1} \}$ and note that $C$ is the $k$-element subset of $X$ with largest sum such that $|A \cap C| = 1$. Using \lref{eigenvecset} we give a lower bound on $b_{C}$, the sum of $C$, under the assumptions that $n \geq k^{2}$ and that there are at most $\binom{n-1}{k-1}$ $k$-element subsets with nonnegative sum in $X$.

\begin{lemma}
\label{lem:packingset}
Let $A = \{x_{1}, \ldots, x_{k} \}$ and let $C = \{x_{1}, x_{k+1}, \ldots, x_{2k-1} \}$. If $n \geq k^{2}$ and there are at most $\binom{n-1}{k-1}$ nonnegative $k$-element subsets of $X$ then $b_{C}$, the sum of $C$, satisfies
\begin{equation}
\label{bcset}
b_{C} \geq \left( 1 - \frac{(2k-1)(k-1)}{n-2k+1} \right) b_{A}.
\end{equation}
\end{lemma}

\noindent \textbf{Proof.} By \lref{eigenvecset} with $j=k-1$, the dot product of the row of $B_{k-1}$ corresponding to $A$ and $\vec{b}$ equals $-(k-1)\binom{n-k}{k-1}b_{A}$. Hence, by  \eqref{bjentries},
\begin{equation}
\label{bk-1}
k \sum_{S \cap A = \emptyset} b_{S} + \sum_{|S \cap A| = 1} b_{S} = -(k-1)\binom{n-k}{k-1}b_{A}.
\end{equation}
Consequently, by \eqref{sumonempty} and Pascal's identity,
\begin{equation}
\label{a1}
\sum_{|S \cap A| = 1} b_{S} = \left( \binom{n-k-1}{k-1} -(k-1)\binom{n-k-1}{k-2} \right) b_{A},
\end{equation}
which is nonnegative if and only if $n \geq k^{2}$. Observe that for any $S \in \binom{X}{k}$ such that $|S \cap A| = 1$, we have $b_{C} \geq b_{S}$. We claim that
\begin{equation}
\label{similartopacking}
b_{C} \geq \left( \frac{ \binom{n-k-1}{k-1} - (k-1)\binom{n-k-1}{k-2} }{ \binom{n-1}{k-1} } \right) b_{A},
\end{equation}
otherwise as $b_{A}$ is a largest entry of $\vec{b}$, \eqref{a1} implies there are at least $\binom{n-1}{k-1}$ nonnegative entries $b_{S}$ where $|S \cap A| = 1$. Since $A$ has nonnegative sum, we would get greater than $\binom{n-1}{k-1}$ nonnegative $k$-element subsets of $X$ if \eqref{similartopacking} does not hold.

Now, we show that the fraction on the right hand of \eqref{similartopacking} is at least the fraction on the right hand side of \eqref{bcset}. We have
\begin{equation}
\label{bcsetsimplify1}
\binom{n-k-1}{k-1} - (k-1) \binom{n-k-1}{k-2} = \left( 1 - \frac{(k-1)^{2}}{n-2k+1} \right) \binom{n-k-1}{k-1}.
\end{equation}
We also have that
\begin{equation}
\label{bcsetsimplify2}
\frac{ \binom{n-k-1}{k-1} }{ \binom{n-1}{k-1} } = \frac{ (n-k-1) \cdots (n-2k+1) }{ (n-1) \cdots (n-k+1) } > \left( 1 - \frac{k}{n-k+1} \right)^{k-1} > 1 - \frac{k(k-1)}{n-k+1}.
\end{equation}
Putting \eqref{bcsetsimplify1} and \eqref{bcsetsimplify2} together yields \eqref{bcset}. \qed

Now we give a lower bound on the number of nonnegative $k$-element subsets that contain $x_{1}$ under the assumptions that $n \geq k^{2}$ and that there are at most $\binom{n-1}{k-1}$ $k$-element subsets of $X$ with nonnegative sum.

\begin{lemma}
\label{lem:lotson1sets}
If $n \geq k^{2}$ and there are at most $\binom{n-1}{k-1}$ nonnegative $k$-element subsets of $X$, then the number of nonnegative $k$-element subsets that contain $x_{1}$ is at least
\begin{equation}
\label{generallotson1set}
\left( 1 - \frac{(6k-3)(k-1)}{n-2k+1} \right) \binom{n-1}{k-1}.
\end{equation}
\end{lemma}

\noindent \textbf{Proof.} Recall that $A = \{x_{1}, \ldots, x_{k} \}$ and that $C = \{x_{1}, x_{k+1}, \ldots, x_{2k-1} \}$. By \lref{eigenvecset}, the dot product of the row of $B_{k}$ corresponding to $C$ and $\vec{b}$ equals $-\binom{n-k-1}{k-1} b_{C}$. Hence,
\begin{equation}
\label{rinseandrepeat}
\sum_{\substack{S \cap C \neq \emptyset, \\ S \cap A \neq \emptyset}} b_{S} + \sum_{\substack{S \cap C \neq \emptyset, \\ S \cap A = \emptyset}} b_{S} = \sum_{\substack{S \cap C \neq \emptyset}} b_{S} = \binom{n-k-1}{k-1} b_{C}.
\end{equation}
We claim that
\begin{equation}
\label{nottoobig}
\sum_{\substack{S \cap C \neq \emptyset, \\ S \cap A = \emptyset}} b_{S} \leq \left( \binom{n-1}{k-1} - \binom{n-k-1}{k-1} \right) b_{A}.
\end{equation}
Otherwise, there would be at least $\binom{n-1}{k-1} - \binom{n-k-1}{k-1}$ nonnegative entries $b_{S}$ such that $S \cap C \neq \emptyset$ and $S \cap A = \emptyset$ as $b_{A}$ is a largest entry. By \lref{highestweightintersectionset}, this would yield greater than $\binom{n-1}{k-1}$ nonnegative $k$-element subsets. Hence, \eqref{nottoobig} holds, which implies by \lref{packingset}, \eqref{bcsetsimplify2}, \eqref{rinseandrepeat}, and \eqref{nottoobig} that
\begin{align}
\label{prettybigset}
\sum_{\substack{S \cap C \neq \emptyset, \\ S \cap A \neq \emptyset}} b_{S} &\geq \left( \left(1 - \frac{k(k-1)}{n-k+1} \right) \left(2 - \frac{(2k-1)(k-1)}{n-2k+1} \right) - 1 \right) \binom{n-1}{k-1} b_{A} \notag \\
& \geq \left( 1 - \frac{(4k-1)(k-1)}{n-2k+1} \right) \binom{n-1}{k-1} b_{A}.
\end{align}

Let $\cF_{i}$ be the family of $k$-element subsets of $X$ that contain $x_{i}$ but not $x_{1}$ and intersect $A$ and $C$. We have
\begin{equation}
\label{whereweregoing}
\sum_{x_{1} \in S} b_{S} = \sum_{\substack{S \cap C \neq \emptyset, \\ S \cap A \neq \emptyset}} b_{S} - \sum_{\substack{S \cap C \neq \emptyset, \\ S \cap A \neq \emptyset, \\ x_{1} \notin S}} b_{S} = \sum_{\substack{S \cap C \neq \emptyset, \\ S \cap A \neq \emptyset}} b_{S} - \sum_{i=2}^{n} |\cF_{i}|x_{i}.
\end{equation}
We first show that if $i \in \{2, \ldots, 2k-1\}$ then
\begin{equation}
\label{case1}
|\cF_{i}| = \binom{n-2}{k-1} - \binom{n-k-1}{k-1}.
\end{equation}
Without loss of generality suppose that $x_{i} \in A \setminus \{x_{1}\}$. There are $\binom{n-2}{k-1}$ $k$-element sets of $X$ that contain $x_{i}$ but not $x_{1}$. From these, we subtract the $\binom{n-k-1}{k-1}$ $k$-element subsets of $X$ that contain $x_{i}$ but do not intersect $C$.

Now we determine $|\cF_{i}|$ when $i \in \{2k, \ldots, n\}$. Let $\mathcal{G}_{i}$ (respectively $\cH_{i}$) be the family of $k$-element subsets of $X$ that contain $x_{i}$ but not $x_{1}$ and intersect $A$ (respectively $C$). We have $\cF_{i} = \cG_{i} \cap \cH_{i}$ so by inclusion-exclusion,
\begin{align}
\label{inclexcl}
| \cF_{i} | &= |\cG_{i} \cap \cH_{i}| = |\cG_{i}| + |\cH_{i}| - |\cG_{i} \cup \cH_{i}| \notag \\
&= 2 \left( \binom{n-2}{k-1} - \binom{n-k-1}{k-1} \right) - \left( \binom{n-2}{k-1} - \binom{n-2k}{k-1} \right) \notag \\
&= \binom{n-2}{k-1} - 2 \binom{n-k-1}{k-1} + \binom{n-2k}{k-1}.
\end{align}

By \eqref{whereweregoing}, \eqref{case1}, and \eqref{inclexcl},
\begin{align}
\label{evenmorebetter}
\sum_{\substack{S \cap C \neq \emptyset, \\ S \cap A \neq \emptyset, \\ x_{1} \notin S}} b_{S} &= |\cF_{2}| \sum_{i=2}^{2k-1} x_{i} + |\cF_{2k}| \sum_{i=2k}^{n} x_{i} = |\cF_{2}|(b_{A} + b_{C} - 2x_{1}) + |\cF_{2k}|(x_{1} - b_{A} - b_{C}) \notag \\
&= (2|\cF_{2}| - |\cF_{2k}|)(-x_{1}) + (|\cF_{2}| - |\cF_{2k}|)(b_{A} + b_{C}) \notag \\
&< 2(|\cF_{2}| - |\cF_{2k}|)b_{A} = 2 \sum_{j=k+2}^{2k} \binom{n-j}{k-2} b_{A} \notag \\
&< 2(k-1) \binom{n-k-2}{k-2}b_{A} < \frac{2(k-1)^{2}}{n-1} \binom{n-1}{k-1} b_{A}.
\end{align}

By \eqref{prettybigset}, \eqref{whereweregoing}, and \eqref{evenmorebetter}, we have
\begin{equation}
\label{almostthereset}
\sum_{x_{1} \in S} b_{S} \geq \left( 1 - \frac{(6k-3)(k-1)}{n-2k+1} \right) \binom{n-1}{k-1} b_{A}.
\end{equation}
Hence, a lower bound on the number of nonnegative $k$-element subsets that contain $x_{1}$ is given by \eqref{generallotson1set}. \qed

\section{Bounds from Averaging}
\label{sec:greenekleitman}

The main result of this section is \lref{nonnegTset}, which shows that if $T \in \binom{X}{k}$ is a $k$-element subset with negative sum, then there are at least $\binom{n-2k}{k-1}$ nonnegative $k$-element subsets of $X$ that are disjoint from $T$. The proof of \lref{nonnegTset} is similar to Manickam's and Singhi's proof of \conjref{MMS} when $k|n$ \cite{ms} and has also been observed by others \cite{ahs, frankl, Tyomkyn}.

\begin{lemma}
\label{lem:nonnegTset}
If $T \in \binom{X}{k}$ has negative sum, then there are at least
\begin{equation}
\label{nonnegtrivialwithTset}
\binom{n-2k}{k-1} \geq \left( 1 - \frac{(2k-1)(k-1)}{n-2k+1} \right) \binom{n-1}{k-1}
\end{equation}
nonnegative $k$-element subsets of $X$ that are disjoint from $T$.
\end{lemma}

\noindent \textbf{Proof.} Write $n = mk + r$ where $0 \leq r \leq k-1$. Since $T \in \binom{X}{k}$ has negative sum, adding the smallest $r$ elements of $X \setminus T$ to $T$ yields a $(k+r)$-element subset $U \in \binom{X}{k+r}$ with negative sum. Let
\begin{equation}
\label{disjointfromU}
\cF := \left \{ S \in \binom{X \setminus U}{k} : b_{S} \geq 0 \right \}
\end{equation}
be the family of $k$-element sets disjoint from $U$ that have nonnegative sum.

Consider a random permutation $\pi \in S_{X}$ that fixes $U$. Partition the $(m-1)k$ elements of $X \setminus U$ into $k$-element sets $\cS = \{S_{1}, \ldots, S_{m-1}\}$, and define the indicator random variable $Z_{i}$ to be $1$ if $\pi(S_{i})$ has nonnegative sum and $0$ otherwise. Let $Z = \sum_{i=1}^{m-1} Z_{i}$ and note that $Z \geq 1$ because some $k$-element subset in $\pi(\cS)$ must have nonnegative sum as the sum of the elements of $X \setminus U$ is positive. On the other hand, $\mathbb{E}(Z_{i})$ is the probability that a randomly chosen $k$-element subset that is disjoint from $U$ has nonnegative sum. Hence,
\begin{equation}
\label{greatexpections}
\mathbb{E}(Z_{i}) = \frac{|\cF|}{\binom{n-k-r}{k}}.
\end{equation}
By linearity of expectation,
\begin{equation}
\label{loe}
1 \leq \mathbb{E}(Z) = (m-1) \mathbb{E}(Z_{i}) = \frac{(n-k-r)|\cF|}{k \binom{n-k-r}{k}}.
\end{equation}
Since $0 \leq r \leq k-1$, we have
\begin{equation}
\label{boundF}
|\cF| \geq \binom{n-k-r-1}{k-1} \geq \binom{n-2k}{k-1} \geq \left( 1 - \frac{(2k-1)(k-1)}{n-2k+1} \right) \binom{n-1}{k-1}.
\end{equation}
Since $U \in \binom{X}{k+r}$ contains $T \in \binom{X}{k}$, each $k$-element subset in $\cF$ is also disjoint from $T$. \qed

\section{Proof of \tref{quadratic}}
\label{sec:thequadraticproof}

Finally, we prove \tref{quadratic}.

\vspace{0.25cm}

\noindent \textbf{Proof of \tref{quadratic}} If all $k$-element subsets containing $x_{1}$ have nonnegative sum, then there are at least $\binom{n-1}{k-1}$ $k$-element subsets of $X$ with nonnegative sum.

Otherwise, some $k$-element subset $T \in \binom{X}{k}$ containing $x_{1}$ has negative sum. Suppose, for a contradiction, that there are at most $\binom{n-1}{k-1}$ nonnegative $k$-element subsets in this case. By \lref{lotson1sets}, there are at least
\begin{equation}
\label{spacesonx1}
\left(1 - \frac{(6k-3)(k-1)}{n-2k+1} \right) \binom{n-1}{k-1}
\end{equation}
nonnegative $k$-element subsets containing $x_{1}$ since $n \geq 8k^{2}$.

Since $T \in \binom{X}{k}$ has negative sum, by \lref{nonnegTset}, there are at least
\begin{equation}
\label{disjointTset}
\left(1 - \frac{(2k-1)(k-1)}{n-2k+1} \right) \binom{n-1}{k-1}
\end{equation}
nonnegative $k$-element subsets of $X$ that have trivial intersection with $T$.

Since $T$ contains $x_{1}$, none of the $k$-element subsets counted in \eqref{disjointTset} contain $x_{1}$. Summing \eqref{spacesonx1} and \eqref{disjointTset}, there are at least
\begin{align}
\label{grandfinaleset}
\left( 2 - \frac{(8k-4)(k-1)}{n-2k+1} \right) \binom{n-1}{k-1}
\end{align}
nonnegative $k$-element subsets in $X$. For $n \geq 8k^{2}$, however, the expression in \eqref{grandfinaleset} is greater than $\binom{n-1}{k-1}$, which contradicts our assumption. \qed

\medbreak
\noindent{\sc Acknowledgement:}  We are grateful to Simeon Ball, Chris Godsil, Po-Shen Loh, Karen Meagher, Bruce Rothschild, Benny Sudakov, Terence Tao, Jacques Verstra\"{e}te, and Rick Wilson for their advice.

\bibliographystyle{plain}
\bibliography{MMSQuadratic}

\end{document}